\documentclass[10pt,a4paper]{amsart}
\usepackage[latin1]{inputenc}
\usepackage{amsmath}
\usepackage{amsfonts}
\usepackage{amssymb}
\usepackage{graphicx}
\usepackage{amscd}
\usepackage{mathrsfs}
\usepackage[all]{xy}
\usepackage{overpic}

\newtheorem{thm}{Theorem}[section]
\newtheorem{lemma}[thm]{Lemma}
\newtheorem{prop}[thm]{Proposition}

\theoremstyle{definition}
\newtheorem{defn}[thm]{Definition}

\theoremstyle{remark}
\newtheorem{remark}[thm]{Remark}
\numberwithin{equation}{section}

\newcommand{\Comp}{\mathbb{C}}       
\newcommand{\Real}{\mathbb{R}}       
\newcommand{\Ham}{\mathbb{H}}        
\newcommand{\R}{\mathbb{R}}      

\begin{document}

    \author{Andrew Clarke}
    \title{Instantons on the Exceptional Holonomy Manifolds of Bryant and Salamon}    
\address{Instituto de Matem\'atica, Universidade Federal do Rio de Janeiro, Av. Athos da Silveira Ramos 149, Rio de Janeiro, RJ, 21941-909, Brazil.}
 \email{andrew@im.ufrj.br}
     
\date{\today}

\maketitle


\begin{abstract}
We give a construction of $G_2$ and $Spin(7)$ instantons on exceptional holonomy manifolds constructed by Bryant and Salamon, by using an ansatz of spherical symmetry coming from the manifolds being the total spaces of rank-$4$ vector bundles. In the $G_2$ case, we show that, in the asymptotically conical model, the connections are asymptotic to Hermitian Yang-Mills connections on the nearly K\"ahler $S^3\times S^3$. 
\end{abstract}


\section{Introduction}

A natural situation in which to study gauge fields is when the Riemannian manifold $(M,g)$ has a special holonomy group, in which case the holonomy group determines a subalgebra $\mathfrak{g}\subseteq\Lambda^2$. For $G_2$, $Spin(7)$ and Calabi-Yau manifolds, if the curvature $F_A$ of a connection lies in $\mathfrak{g}$, the connection automatically satsfies the Yang-Mills equations. 
 
Following the constructions of Bryant and Salamon, and of Joyce, the study of gauge fields on reduced holonomy manifolds was formulated in the physics literature (see for example \cite{AOS,BKS}) and also in mathematics (see \cite{DT,Ti}). Existence theorems have been given on both compact and non-compact spaces in, for example, \cite{S,W,FN,HILP,IP}. 
 

Here we consider the problem of constructing solutions to the instanton equation on the non-flat manifolds of reduced holonomy that were obtained by Bryant and Salamon. In \cite{BS} a number of  complete metrics of holonomy exactly $G_2$ and $Spin(7)$ were constructed. Their work stemmed from a principal bundle construction to obtain a family of non-degenerate differential forms  that depend on $2$  functions of one variable. In each given case they were able to solve a system of ordinary differential equations so as to obtain a metric $g_\gamma$ of exceptional holonomy. 

The differentiable manifolds that support the metrics are the total spaces of rank $3$ and $4$ vector bundles over $3$ and $4$-dimensional manifolds.
In one example for the group $G_2$, the manifold $X$ is given as the total space of the spinor bundle $ \mathcal{S}\to S^3$ over the three-sphere. The metric $g_\gamma$ is then given by
\begin{eqnarray*}
g_\gamma = 4(1+r)^{-1/3}d\sigma^2 +3(1+r)^{2/3}ds^2
\end{eqnarray*}
where $ds^2$ is the pullback of the round metric tensor on $S^3$ and $d\sigma^2$ is a quadratic form on $X$ that restricts to the fibres of the spinor bundle to be the standard flat metric. That is, $g_\gamma$ restricts to the fibres to be conformally equivalent to the flat metric on $\R^4$.

The original ansatz used to obtain instanton solutions to the Yang-Mills equations in $4$-dimensions was to use the assumption of spherical symmetry, and to construct 
\begin{eqnarray}\label{eqn:BPSTinstanton}
A=\frac{\text{Im}(xd\bar x)}{1+|x|^2}
\end{eqnarray}
so that $\nabla=d+A$ is an $SU(2)$-connection with self-dual curvature on $\R^4$ (see \cite{BPST}). 

In the example at hand, we take the $\mathfrak{su}(2)$-valued form 
\begin{eqnarray*}
A_1=\text{Im}(a\bar\alpha)
\end{eqnarray*}
 and construct a connection whose curvature lies in the subspace $\mathfrak{g}_2\otimes\mathfrak{su}(2)\subseteq \Lambda^2\otimes \mathfrak{su}(2)$. Neither the quaternion-valued function $a$ or form $\alpha$ are well defined on $X$, but the form $A_1$ is well-defined. We set $r=|a|^2=a\bar a$ to be the squared fibre radius function on $\mathcal{S}$. This function is well-defined on $X$. 

\begin{thm}\label{thm:G2}
Let $X$ be the total space of the spinor bundle over the round $S^3$ of constant curvature $\kappa$, equipped with the $G_2$-holonomy metric that is determined by the form $\gamma$ constructed by Bryant and Salamon. For any $C>-3$ let $f(r)$ be the function
\begin{eqnarray*}
f(r)=\frac{2}{3(r+1)+C(r+1)^{1/3}}.
\end{eqnarray*}
Then the connection $\nabla=d+f(r)A_1$ defines a non-trivial $G_2$-instanton on the trivial rank $2$ complex vector bundle over $X$. The connections are invariant under the group $SU(2)^3$ of isometries of $X$.
\end{thm}
The first statements of this result are proven in Section \ref{sec:G2inst}. The statement on the invariance is shown in Section \ref{sec:G2asymp}.

It has also been noted that the metrics of Bryant and Salamon are asymptotically conical. That is, outside of sufficiently large compact sets, the metric $g_\gamma$ is arbitrarily close to a conical metric
\begin{eqnarray*}
d\rho^2 +\rho^2 g 
\end{eqnarray*}
 defined on $\R^+\times (S^3\times S^3)$. The closeness is measured with respect to the conical metric. A consequence of $(X,g_\gamma)$ having $G_2$ holonomy is that $(S^3\times S^3,g)$ must be \emph{nearly K\"ahler}. This is to say that it admits an almost-complex structure compatible with the metric and a complex $(3,0)$-form that satisfies certain differential relations with the Hermitian form of the metric. In this context we can consider the Hermitian Yang-Mills equations for connections on principal and vector bundles over $S^3\times S^3$. We have the following result on the asymptotic behaviour of the connections constructed in Theorem \ref{thm:G2}.
\begin{thm}\label{thm:asymptotics}
There exists a non-flat Hermitian Yang-Mills connection $\widetilde{A}$ on the trivial $\Comp^2$ bundle over $S^3\times S^3$ such that every $G_2$-instanton on $X$ constructed in Theorem \ref{thm:G2} is asymptotic to the pull-back of $\widetilde{A}$ to $\Real^+\times (S^3\times S^3)$.
\end{thm}
This result is demonstrated in Section \ref{sec:G2asymp}.

To the author's knowledge, this example was previously unknown, though it would be interesting to understand it in comparison with the works of Bryant, and Harland and N\"olle \cite{B,HN}.

The results in this paper have much in common with the recent work of Oliveira \cite{O}. The principal difference is that he considers monopoles and instantons on the other $G_2$ metrics of Bryant and Salamon, those being total spaces of rank-3 vector bundles over $4$-manifolds. The convenient model for the metric in that case is that of \cite{CGLP}, rather than the original model of Bryant and Salamon. 

The Bryant-Salamon construction for metrics of $Spin(7)$-holonomy involves considering the negative spinor bundle on $S^4$ and constructing a non-degenerate $4$-form. Bryant and Salamon solved a system of ODE's to obtain an example for which the form $\Psi$ is closed. We are able to construct a $Spin(7)$-instanton with structure group $SU(2)$ for this $Spin(7)$ structure. 

Let $\mathcal{S}^-\to S^4$ be the negative spinor bundle on $S^4$ and denote by $Y^8$ the total space of this manifold. From the bundle $\mathcal{S}^-$, we can construct a principal $SU(2)$-bundle $\mathcal{E}_Y$ over $Y$ and consider connection forms on $\mathcal{E}_Y$ of the form
\begin{eqnarray*}
A=\phi+f(r)A_2
\end{eqnarray*}
where $\phi$ is pulled back from the connection on $\mathcal{S}^-$ for the round metric on $S^4$, and $A_2=\text{Im}(\bar a\alpha)$ is an equivariant $1$-form on $\mathcal{E}_Y$. $f(r)$ is a function on $X$ that depends only on the radius in the $\mathcal{S}^-$-fibre directions. 
\begin{thm}\label{thm:Spin7}
Let $Y$ be the total space of the negative spinor bundle on $S^4$, equipped with the $Spin(7)$-holonomy metric of Bryant and Salamon.  Let $\mathcal{E}_Y$ be the associated $SU(2)$-bundle on $Y$ and let $f(r)$ be the function  
\begin{eqnarray*}
f(r)&=& \frac{1}{r(1+D(1+r)^{3/5})}+\frac{D(2r+5)}{5r(1+r)^{2/5}(1+D(1+r)^{3/5})}
\end{eqnarray*}
where $D>-1$ is a constant. Then the connection $A=\phi+fA_2$ defines a non-trivial $Spin(7)$-instanton on $\mathcal{E}_Y$ with structure group $SU(2)$. 
\end{thm}
 The connection form $A$ is defined only away from the set $\{r=0\}$ in $Y$, that being the image of the zero section of $\mathcal{S}^-$. This theorem is proved in Section \ref{sec:Spin7inst}.

 
\section{Gauge theory on $G_2$ manifolds} 

We now consider manifolds of $G_2$ holonomy and derive the first order instanton equation in this context. 

For any Riemannian metric $g$, the condition that $Hol(g)\subseteq G_2$ is equivalent to the existence of a non-degenerate $3$-form $\varphi$ such that 
\begin{eqnarray} \label{eqn:phi parallel}
\nabla^g\varphi =0.
\end{eqnarray}
Any non-degenerate $3$-form $\varphi$ algebraically determines a unique metric $g_\varphi$ and if (\ref{eqn:phi parallel}) holds then the uniqueness of the Levi-Civita connection implies that $\nabla^g=\nabla^\varphi$ is the Levi-Civita connection for the metric $g_\varphi$. Then, the non-linear equation in $\nabla^\varphi\varphi=0$ is equivalent to the equations
\begin{eqnarray*}
d\varphi &=&0\\
d(*_\varphi\varphi)&=&0
\end{eqnarray*}
where $*_\varphi$ is the Hodge dual operator for the metric $g_\varphi$.

 The reducibility of the action of $G_2$ on $\Lambda^2\R^7$ implies that on a $G_2$-manifold $X$ there exists a decomposition 
\begin{eqnarray}\label{Lam2decomposition}
\Lambda^2T^*_X=\Lambda^2_7\oplus \Lambda^2_{14}
\end{eqnarray}
into the direct sum of sub-bundles of ranks $7$ and $14$ respectively. $\Lambda^2_{14}$ is equal to the kernel of the map
\begin{eqnarray*}
*\varphi\wedge \cdot:\Lambda^2\to \Lambda^6.
\end{eqnarray*}
The subspaces $\Lambda^2_7$ and $\Lambda^2_{14}$ are eigenspaces of the map $\alpha\mapsto *(\varphi\wedge\alpha)$. With the conventions of \cite{BS},
\begin{eqnarray*}
\alpha\in\Lambda^2_7 &\iff & *(\varphi\wedge \alpha)=-2\alpha\\
\alpha\in\Lambda^2_{14} &\iff & *(\varphi \wedge \alpha)=\alpha.
\end{eqnarray*}

Similarly, the holonomy of an $8$-dimensional Riemannian manifold $(Y,g)$ being contained in $Spin(7)$ is equivalent to the existence of a non-degenerate $4$-form $\Psi$ that satisfies $d\Psi=0$. $Spin(7)$ acts reducibly on $\Lambda^2\R^8=\Lambda^2_7\oplus\Lambda^2_{21}$ and each summand is an eigenspace for a linear map. That is,
\begin{eqnarray*}
\beta\in\Lambda^2_7 &\iff& *(\Psi\wedge\beta)=-3\beta\\
\beta\in\Lambda^2_{21} &\iff& *(\Psi\wedge\beta)=\beta
\end{eqnarray*}
where $*$ is the Hodge operator $*:\Lambda^6\to \Lambda^2$.

\begin{defn}
Let $E$ be a vector bundle over a $G_2$-manifold $(X,\varphi)$. A $G_2$-instanton on $E$ is a connection $A$ whose curvature satisfies  $*\varphi\wedge F_A=0$.
\end{defn}

\begin{remark}
The factors in the decomposition (\ref{Lam2decomposition}) are each eigenspaces (at each point) for the operator $\alpha\mapsto *(\varphi\wedge\alpha)$ so one could equally well define an instanton to be a connection $A$ that satisfies $*(\varphi\wedge F_A)=F_A$.
\end{remark}

\begin{defn} 
Let $E$ be a vector bundle over a $Spin(7)$-manifold $(Y,\Psi)$. A $Spin(7)$-instanton on $E$ is a connection $A$ whose curvature satisfies $*(\Psi\wedge F_A)=F_A$.
\end{defn}

 
\section{The manifolds of Bryant and Salamon}    \label{Sect:BSmetric}

Let $(M,g)$ be a $3$-dimensional Riemannian manifold, with $\mathcal{F}\to M$ the principal coframe bundle. $\mathcal{F}$ is then equipped with a canonical $\text{Im}\Ham\cong\mathfrak{su}(2)$-valued $1$-form $\omega$ and Levi-Civita connection form $\phi$. Let $\tilde{\mathcal{F}}$ be a spin structure for $M$, in the sense that $\tilde{\mathcal{F}}$ is a double cover of $\mathcal{F}$ compatible with the covering of $SO(3)$ by $SU(2)$. The form $\omega$ on $\tilde{\mathcal{F}}$ satisfies $R_g^*\omega=\bar g\omega g=Ad_{g^{-1}}\cdot \omega$.

Assume that $M$ has constant positive sectional curvature. This means that $\omega$ and $\phi$ satisfy the relation 
\begin{eqnarray*}
d\phi+\phi\wedge \phi =-\frac{\kappa}{2}\Omega
\end{eqnarray*}
where $\Omega=\frac{1}{2}\omega\wedge\omega = \omega^{23}i+\omega^{31}j+\omega^{12}k$ and where $\kappa>0$ is a constant. 
The bundle of spinors on $M$ is given as $\mathcal{S}=(\tilde{\mathcal{F}}\times\Ham)/SU(2)$ where $SU(2)$ acts on the two factors on the right. Let $X=\mathcal{S}$ be the $7$-dimensional manifold given as the total space of this bundle. Then $\mathcal{E}_X=\tilde{\mathcal{F}}\times \Ham$ can be considered a principal $SU(2)$-bundle over $X$.

The forms $\omega$ and $\phi$ can be considered forms on $\mathcal{E}_X$ with values in $\text{Im}\Ham$. Also, we consider the map $a:\tilde{\mathcal{F}}\times\Ham\to \Ham$ given by projection to the second factor and the form $\alpha=da-a\phi$ on $\mathcal{E}_X$. 

Then we write $\omega =\omega^1 i+\omega^2 j+\omega^3k$ and $\alpha=\alpha^0+\alpha^1i+\alpha^2j+\alpha^3k$ and 
\begin{eqnarray*}
B=1/2\bar{\alpha}\wedge\alpha =B^1i+B^2j+B^3k =(\alpha^{01}-\alpha^{23})i+(\alpha^{02}-\alpha^{31})j+(\alpha^{03}-\alpha^{12})k.
\end{eqnarray*}
The form $\alpha$ satisfies $R_g^*\alpha=\alpha g$ so  we have $R_g^* B=1/2R_g^*(\bar\alpha\wedge \alpha)= \bar g\cdot B\cdot g =Ad_{g^{-1}}\cdot B$. We can define $3$ and $4$-forms on $\mathcal{E}_X$
\begin{eqnarray*}
\gamma_1&=&\omega^1\wedge \omega^2\wedge\omega^3=\frac{-1}{6}\omega\wedge\omega\wedge\omega\\
\gamma_2&=& \omega^1\wedge B^1 +\omega^2 \wedge B^2 +\omega^3\wedge B^3\\
&=& - \text{Re}(\omega\wedge B)\\
\psi_1 &=& \alpha^0\wedge \alpha^1\wedge \alpha^2\wedge \alpha^3\\
\psi_2&=& \omega^{23}\wedge B^1 + \omega^{31}\wedge B^2 +\omega^{12}\wedge B^3 \\
&=&-\text{Re}(\Omega\wedge B)
\end{eqnarray*}
The forms $\gamma_1$, $\gamma_2$, $\psi_1$, $\psi_2$ vanish when contracted by elements of the kernel of $\pi_*:T\mathcal{E}_X\to TX$ and are invariant under the action of $SU(2)$ on $\mathcal{E}_X$. For example, $\gamma_2=\langle\omega, B\rangle$ is given by the $Ad$-invariant Killing form on $\mathfrak{su}(2)$, evaluated on the forms $\omega$ and $B$. Then,
\begin{eqnarray*}
R_g^*\gamma_2=\langle Ad_{g^{-1}}\omega, Ad_{g^{-1}} B\rangle =\langle\omega,B\rangle= \gamma_2.
\end{eqnarray*}
 As such they are each pull-backs to $\mathcal{E}_X$ of $3$ and $4$-forms on $X$.


Then, for smooth functions $f$ and $g$, depending only upon the squared radial function $r=|a|^2$, we set $\gamma=f^3\gamma_1 +fg^2\gamma_2$. This is a \emph{non-degenerate} $3$-form in the sense of $G_2$-geometry and so determines a  metric and orientation with respect to which the Hodge dual of $\gamma$ is given by
\begin{eqnarray*}
*_\gamma \gamma =g^4\psi_1-f^2g^2\psi_2.
\end{eqnarray*}
Then,
\begin{eqnarray*}
d\gamma &=&\left[ (f^3)^\prime -(\frac{3\kappa}{4} )fg^2\right]  dr\wedge \gamma_1 +(fg^2)^\prime dr\wedge \gamma_2\\
d*_\gamma\gamma &=&\left[ -(f^2g^2)^\prime +(\frac{\kappa}{4})g^4\right] dr\wedge \psi_2
\end{eqnarray*}
Bryant and Salamon show that for the functions
\begin{eqnarray*}
f(r)&=&(3\kappa)^{1/2}(1+r)^{1/3}\\
g(r)&=& 2(1+r)^{-1/6}
\end{eqnarray*}
the form $\gamma$ is closed and coclosed and hence determines a torsion-free $G_2$ structure on $X$. Furthermore, Bryant and Salamon show that the metric is complete and irreducible. In particular, 
\begin{eqnarray*}
\psi=*_\gamma\gamma &=& \left(\frac{16}{(1+r)^{2/3}}\right)\psi_1+
\left(-12\kappa(1+r)^{1/3}\right)\psi_2\\
&=&\sigma\psi_1+\tau\psi_2.
\end{eqnarray*}
For future reference, we note that
\begin{eqnarray}
\frac{-\kappa\sigma}{4\tau}=\frac{-16\kappa(1+r)^{-2/3}}{-48\kappa(1+r)^{1/3}}=\frac{1}{3}\frac{1}{1+r}.
\end{eqnarray}



In our second example, we consider the complete $8$-dimensional Riemannian manifold with $Spin(7)$ holonomy that is constructed by Bryant and Salamon. The manifold $Y$ is the total space of the negative spinor bundle over the round $4$-sphere. That is, just as in the previous example, $Y$ is the total space of a rank-$4$ vector bundle over a space-form.

Let $\pi:\tilde{\mathcal{F}}\to S^4$ be the principal coframe bundle over $M=S^4$ and let $\mathcal{F}$ be the spin double cover of $\tilde{\mathcal{F}}$. That is, $\mathcal{F}$ is a principal bundle over $S^4$ with structure group $Spin(4)=SU(2)\times SU(2)$. $\mathcal{F}$ canonically admits an $\Ham$-valued $1$-form $\omega$ and we set $\Omega=\frac{1}{2}\bar\omega\wedge\omega$. 
The Levi-Civita connection on $\tilde{\mathcal{F}}$ consists of two $\mathfrak{su}(2)$-valued $1$-forms $\xi$ and $\phi$ that satisfy 
\begin{eqnarray*}
d\omega=-\xi\wedge\omega-\omega\wedge\phi.
\end{eqnarray*}
and the Riemannian curvature tensor is given in the form of $\mathfrak{su}(2)$-valued forms $\Phi=d\phi+\phi\wedge\phi$ and $\Xi=d\xi+\xi\wedge\xi$. The condition that $S^4$ is Einstein and with self-dual Weyl curvature is equivalent to the identity
\begin{eqnarray*}
d\phi+\phi\wedge\phi =\frac{\kappa}{2}\Omega
\end{eqnarray*}
where the constant $\kappa>0$ is a multiple of the scalar curvature. 

 We consider the representation of $Spin(4)$ on $\Ham$ given by
\begin{eqnarray*}
 ((p,q), v)\mapsto v\bar q
\end{eqnarray*}
The associated vector bundle $\mathcal{S}^-=(\mathcal{F}\times \Ham)/Spin(4)$ is the negative spinor bundle of $S^4$. We denote this space by $Y$ when considering this as a smooth $8$-dimensional manifold. Let $a:\mathcal{F}\times \Ham\to \Ham$ denote the projection onto $\Ham$ and let $\alpha=da-a\phi$. The basic $1$-forms on $\mathcal{F}\times \Ham$, with respect to the projection to $Y$ are then spanned at each point by the components of the forms $\omega$ and $\alpha$. 

We consider the differential forms on $\mathcal{F}\times \Ham$ :
\begin{eqnarray*}
\Omega=\frac{1}{2}\bar\omega\wedge\omega 
&=& i(\omega^{01}-\omega^{23})+j(\omega^{02}-\omega^{31})+k(\omega^{03}-\omega^{12})\\
&=& i\Omega^1+j\Omega^2+k\Omega^3\\
\frac{1}{2}\bar\alpha\wedge\alpha &=& B =iB^1+jB^2+kB^3\\
\Psi_1 &=& \frac{-1}{6}B\wedge B =\alpha^{0123}\\
\Psi_2 &=& -\text{Re}(B\wedge \Omega) = B^1\Omega^1+B^2\Omega^2+B^3\Omega^3\\
\Psi_3 &=& \frac{-1}{6}\Omega\wedge\Omega =\omega^{0123}.
\end{eqnarray*}
As in the $G_2$ case, it can be noted that the forms $\Psi_1,\Psi_2,\Psi_3$ are invariant under the action of $Spin(4)$ and so descend to forms on $Y$. For any positive functions $\sigma,\tau$ on $Y$, the form 
\begin{eqnarray*}
\Psi=\sigma^2\Psi_1+\sigma\tau\Psi_2+\tau^2\Psi_3
\end{eqnarray*}

is \emph{non-degenerate} in the sense that at every point it can be identified with the fundamental $4$-form on $\R^8$. $\Psi$ defines a reduction of the structure group of $X$ to $Spin(7)$ and the associated metric is given by 
\begin{eqnarray*}
ds^2 =\sigma\big((\alpha^0)^2+(\alpha^1)^2+(\alpha^2)^2+(\alpha^3)^2\big)+\tau\big((\omega^0)^2+(\omega^1)^2+ (\omega^2)^2+(\omega^3)^2\big).
\end{eqnarray*}
As in the $G_2$ case, Bryant and Salamon took the ansatz that $\sigma$ and $\tau$ depend only on $r=a\bar a$. They were able to show that for
\begin{eqnarray*}
\sigma(r)&=& 4(1+r)^{-2/5},\\
\tau(r) &=& 5\kappa (1+r)^{3/5}
\end{eqnarray*}
the form $\Psi$ is closed, and so the holonomy is contained in $Spin(7)$. They furthermore show that the holonomy of this metric is equal to $Spin(7)$. 

$\mathcal{F}\times \Ham$ is a principal $Spin(4)$-bundle over $Y$ and $\mathcal{E}_Y=(\mathcal{F}\times \Ham/(SU(2)\times \{1\})$ is a principal $SU(2)$-bundle over $Y$. The form $\phi$ is invariant under the action of $SU(2)\times \{1\}$ and so descends to define a connection form on $\mathcal{E}_Y$ with curvature given by $\Phi=\kappa/2\Omega$. We will use $\phi$ as a base-point in the set of connections on $\mathcal{E}_Y$ to construct a $Spin(7)$-instanton.


\section{$G_2$-Instantons on the spinor bundle of $S^3$}\label{sec:G2inst}

The manifold $X$ constructed by Bryant and Salamon supports a principal fibre bundle $\mathcal{E}_X$. In this section we construct a connection $A$ on $\mathcal{E}_X$ whose curvature satisfies $*\gamma\wedge F_A=0$. That is, $A$ is a $G_2$-instanton. The curvature is given by 
\begin{eqnarray*}
F_A=dA+A\wedge A.
\end{eqnarray*}
 Let
\begin{eqnarray*}
A_1&=&\text{Im}(a\bar{\alpha})\\
A_2 &=& a\omega\bar{a}.
\end{eqnarray*}
These $1$-forms, which take values in $\mathfrak{su}(2)$, are invariant under the action of $SU(2)$ on $\mathcal{E}_X$ and vanish on the kernel of the projection $\pi_*:T\mathcal{E}_X\to TX$. They are therefore pull-backs of forms on $X$. We consider the connection $A=fA_1+gA_2$ where $f$ and $g$ are functions only of the squared radius function $r=a\bar a$. We then have
\begin{eqnarray*}
F_A&=& f^\prime dr\wedge A_1+fdA_1 +g^\prime dr\wedge A_2+gdA_2\\
&&\ \ \ \ +f^2 A_1\wedge A_1 +fg(A_1\wedge A_2+A_2\wedge A_1)+g^2 A_2\wedge A_2.
\end{eqnarray*}
These forms are well behaved when we differentiate and multiply them. 
\begin{lemma}
\begin{eqnarray*}
dA_1&=& \alpha\wedge\bar\alpha -\frac{\kappa}{2}a\Omega\bar a,\\
dA_2 &=&\alpha\wedge\omega\bar a-a\omega\wedge\bar\alpha=2\text{Im}(\alpha\wedge\omega\bar a) ,\\
A_1\wedge A_1 &=& \frac{-r}{2}\alpha\wedge\bar\alpha -a(\frac{1}{2}\bar \alpha\wedge\alpha)\bar a,\\
A_2\wedge A_2 &=& 2ra\Omega\bar a,\\
A_1\wedge A_2 +A_2\wedge A_1 &=& dr\wedge A_2 -rdA_2,\\
dr\wedge A_1 &=& \frac{r}{2}\alpha\wedge\bar\alpha -a\frac{\bar\alpha\wedge\alpha}{2}\bar a.
\end{eqnarray*}
\end{lemma}


\proof{The proofs of these assertions follows very quickly from the definitions. The only delicate point stems from the non-commutation of forms and functions when they are quaternion-valued. In the case at hand this can be dealt with using the form
\begin{eqnarray*}
dr&=& \alpha\bar a+a\bar \alpha\\
&=& \bar \alpha a+\bar a\alpha
\end{eqnarray*}
and noting that $dr$ is a real-valued $1$-form so anti-commutes with all $\text{Im}\Ham$-valued $1$-forms. 
For example, we have
\begin{eqnarray*}
2dA_1 &=& da\bar \alpha +ad\bar\alpha -d\alpha\bar a +\alpha d\bar a\\
&=&\left[\alpha +a\phi\right]\bar \alpha+a\left[-\phi\bar \alpha -\frac{\kappa}{2}\Omega\bar a\right]-\left[-\alpha\phi+\frac{\kappa}{2}a\Omega\right]\bar a +\alpha\left[\bar \alpha-\phi\bar a\right]\\
&=& 2\alpha\wedge\bar\alpha-\kappa a\Omega\bar a.
\end{eqnarray*}
We next note that if $\beta=\text{Im}\varphi$ is an $\text{Im}\Ham$-valued $1$-form, then $\beta\wedge\beta=\varphi\wedge\varphi=\bar\varphi\wedge\bar\varphi$. We then have
\begin{eqnarray*}
4A_1\wedge A_1 &=& (a\bar\alpha-\alpha\bar a)\wedge(a\bar\alpha-\alpha\bar a)\\
&=& 2A_1\wedge A_1 -r\alpha\wedge\bar\alpha -a\bar\alpha\wedge\alpha\bar a
\end{eqnarray*}
by taking $\varphi=a\bar\alpha$.
}

 The curvature of $A$ is then equal to
\begin{eqnarray}
F_A&=&\left[rf^\prime +2f -rf^2\right] \frac{\alpha\wedge\bar\alpha}{2} + \left[-f^\prime -f^2\right] \frac{a\bar\alpha\wedge\alpha\bar a}{2}\label{eqn:G2curvature}\\
&& \ \ +\left[\frac{-\kappa}{2}f+2rg^2\right]a\Omega\bar a+\left[g^\prime +fg\right]dr\wedge A_2\nonumber\\
&&\ \ \ \ \ \ +\left[g-fgr\right]dA_2.\nonumber
\end{eqnarray}
We denote by $F_1,\ldots,F_5$ the terms in this expression, omitting the functional multiples. That is,
\begin{eqnarray*}
F_1 &=&\frac{1}{2}\alpha\wedge\bar\alpha\\
F_2&=& \frac{1}{2}a\bar\alpha\wedge\alpha\bar a\\
F_3 &=& a\Omega\bar a\\
F_4 &=& dr\wedge A_2 \\
F_5 &=& dA_2=2Im(\alpha\wedge\omega\bar a).
\end{eqnarray*}
The first equation is simplified in this form because we see that $\varphi\wedge\varphi=\bar\varphi\wedge\bar\varphi$ when $\varphi$ is a quaternion-valued $1$-form. 

One can see that these forms are pointwise linearly independent as follows.  By considering self-duality and anti-self-duality in the $\alpha$ variable, and the number of directions in the $\alpha$ and $\omega$ variables respectively it is clear that the only possible linear relation between the $F_i$ terms could be between $F_4$ and $F_5$. At the point $a=1$, $F_4=dr\wedge A_2=2\alpha^0\wedge\omega$ and 
\begin{eqnarray*}
F_5=2Im(\alpha\wedge\bar a)= 2\alpha^0\wedge \omega +2i(\alpha^2\omega^3-\alpha^3\omega^2)+2j(\alpha^3\omega^1-\alpha^1\omega^3)+ 2k(\alpha^1\omega^2-\alpha^2\omega^1)
\end{eqnarray*}
and one sees that the forms are independent for $a=1$. From their evident invariance properties, one sees that they are similarly linearly independent for all other values of $a$. 


We recall the differential form $\psi=\sigma\psi_1+\tau\psi_2=*_\gamma\gamma$ that is constructed by Bryant and Salamon, where $\sigma$ and $\tau$ are functions. The forms $\psi_1$ and $\psi_2$ are given by
\begin{eqnarray*}
\psi_1&=&\alpha^0\wedge\alpha^1\wedge\alpha^2\wedge\alpha^3,\\
\psi_2&=& \omega^{23}\wedge B^1+\omega^{31}\wedge B^2+\omega^{12}\wedge B^3\\
&=& -\text{Re}\left(\Omega\wedge B\right).
\end{eqnarray*}
 The connection $A$ is a $G_2$-instanton if $F_A\wedge \psi=0$. 


\begin{lemma} $\psi_1\wedge F_3= a\Omega\bar a\wedge \alpha^{0123}$. For the other values of $i$, $\psi_1\wedge F_i=0$.
\end{lemma}

\begin{lemma} The curvature terms satisfy

\begin{enumerate}
\item $\psi_2\wedge F_1=0$,
\item $\psi_2\wedge F_2=-2a\Omega\bar a\wedge \alpha^{0123}$,
\item $\psi_2\wedge F_3=0$,
\item $\psi_2 \wedge F_4 = \frac{r}{3}\text{Im}(\alpha\wedge\bar\alpha\wedge\alpha\bar a)\wedge \omega^{123}$,
\item $\psi_2\wedge F_5 = \psi_2\wedge dA_2 = \text{Im}(\alpha\wedge\bar\alpha\wedge\alpha \bar a)\wedge\omega^{123}$.
\end{enumerate}
\end{lemma}

\proof{ It is clear that $\psi_2\wedge F_1=\psi_2\wedge F_3=0$. Also,
\begin{eqnarray*}
\psi_2\wedge F_2&=& (\omega^{23}\wedge B^1+\omega^{31}\wedge B^2+\omega^{12}\wedge B^3)\wedge a\left(\frac{1}{2}\bar \alpha\wedge\alpha\right)\bar a\\
&=& -2a(\omega^{23}i+\omega^{31}j+\omega^{12}k)\bar a\wedge\alpha^{0123}.
\end{eqnarray*}
We have that $A_2=a\omega\bar a$ so 
\begin{eqnarray*}
\psi_2\wedge A_2&=& (\omega^{23}\wedge B^1+\omega^{31}\wedge B^2+\omega^{12}\wedge B^3)\wedge a\omega\bar a\\
&=& a(B^1i+B^2j+B^3k)\bar a\wedge \omega^{123}\\
&=& \frac{1}{2}a\bar\alpha\wedge\alpha\bar a\wedge \omega^{123}\\
\text{so}\ \ \ \psi_2\wedge F_4 &=& \psi_2\wedge dr\wedge A_2\\
&=&\frac{-1}{2}a\bar\alpha\wedge dr\wedge \alpha\bar a\wedge\omega^{123}\\
&=& \frac{-1}{2} a\bar \alpha\wedge (a\bar\alpha +\alpha\bar a)\wedge\alpha\bar a\wedge \omega^{123}\\
&=&-\text{Im}(a\bar\alpha\wedge a\bar\alpha\wedge\alpha\bar a)\wedge\omega^{123}\\
&=& \frac{r}{3}\text{Im}(\alpha\bar\alpha\alpha\bar a)\wedge\omega^{123}.
\end{eqnarray*}
Here we have used the identity $\text{Im}(\bar\varphi\varphi\bar\varphi)=-3\text{Im}(\varphi\varphi\bar\varphi)$ for any $\Ham$-valued $1$-form, which can be quickly verified. 
Finally, 
\begin{eqnarray*}
\psi_2\wedge F_5= \psi_2\wedge dA_2&=& \psi_2\wedge (\alpha\omega\bar a-a\omega\bar \alpha)\\
&=& \alpha(B^1i+B^2j+B^3k)\bar a\wedge \omega^{123}-a(B^1i+B^2j+B^3k)\wedge\omega^{123}\wedge\bar\alpha\\
&=& \frac{1}{2}\alpha\bar\alpha\alpha \bar a\omega^{123}+\frac{1}{2}a\bar\alpha\alpha\bar\alpha\omega^{123}\\
&=& \text{Im}(\alpha\bar\alpha\alpha \bar a)\wedge\omega^{123}.
\end{eqnarray*}

 The above calculations have proved the following lemma.

\begin{lemma}
The $\mathfrak{su}(2)$-valued differential $6$-form $\psi\wedge F_A$ satisfies 
\begin{eqnarray*}
\psi\wedge F_A&=& (\sigma\psi_1+\tau\psi_2)\wedge F_A\\
&=&  \left[\sigma(g^2r-\frac{f\kappa}{4})+\tau(f^\prime +f^2)\right]\Phi_1\\
 &&\ \ + \tau\left[ (g^\prime +fg)\frac{r}{3}+(g-rfg)\right]\Phi_2
\end{eqnarray*}
for $\Phi_1=2a\Omega\bar a\wedge\alpha^{0123}$ and $\Phi_2=\text{Im}(\alpha\bar\alpha\alpha\bar a)\wedge\omega^{123}$.
In particular, $A$ is a $G_2$-instanton if $f$ and $g$ satisfy the system
\begin{eqnarray*}
f^\prime +f^2+\frac{\sigma}{\tau}(g^2r-\frac{f\kappa}{4})&=&0\\
g^\prime +fg+\frac{3}{r}(g-rfg)&=&0.
\end{eqnarray*}
\end{lemma}

This system potentially gives a large family of solutions to the instanton equation. The complication in developing this idea is that this non-linear system is difficult to solve explicitly, or even to determine the asymptotic behavior of any solutions. We simplify the situation by observing that specifying $g\equiv 0$ satisfies the second equation. That is, a connection $A=fA_1$
is a $G_2$-instanton if $f$ satisfies the ordinary differential equation
\begin{eqnarray*}
f^\prime +f^2-\frac{\kappa\sigma}{4\tau}f=0.
\end{eqnarray*}
We recall that in the case of the manifold of $G_2$-holonomy of Bryant and Salamon,
\begin{eqnarray*}
\frac{-\kappa\sigma}{4\tau}=\frac{1}{3}\frac{1}{1+r}
\end{eqnarray*}
so $f$ must satisfy 
\begin{eqnarray}
f^\prime +f^2+\frac{1}{3}\frac{1}{r+1}f=0.\label{eqn:finalnonlinODE}
\end{eqnarray}
This is a Riccati-type equation, and can hence be readily studied. If we pre-suppose that $f=y^\prime/y$ then $y$ satisfies the equation
\begin{eqnarray*}
\frac{1}{y}\left(y^{\prime\prime}+\frac{1}{3}\frac{1}{r+1}y^\prime\right)=0
\end{eqnarray*}
which has solution
\begin{eqnarray*}
y^\prime &=&C_1(r+1)^{-1/3},\\
y&=&3/2C_1(r+1)^{2/3}+C_2.
\end{eqnarray*}
This leads to the general solution to (\ref{eqn:finalnonlinODE})
\begin{eqnarray*}
f(r)=\frac{2}{3(r+1) +C(r+1)^{1/3}}.
\end{eqnarray*}
If the constant $C$ is greater than $-3$, this function is smooth for all (non-negative)-values of $r$, and the connection $A$ is defined on all of $X$. \\


\section{Asymptotic behaviour of the $G_2$-instantons}\label{sec:G2asymp}

In this section we use $r$ to denote the radius function, so that $r^2=a\bar a$. This changes the formulae for the metric and connection, by replacing $r$ by $r^2$. In the following section, we will revert to the previous choice, and to the notation of Bryant and Salamon. Also, for simplicity we consider only the case of the constant $\kappa=1$ in this section. 

The metrics of Bryant and Salamon were also discovered, almost contemporaneously, by Gibbons, Page and Pope \cite{GPP} and the $G_2$-metrics are also studied in detail by Atiyah and Witten \cite{AW}. In those papers the asymptotic behaviour of the metrics is emphasised, and it is shown that the metrics are asymptotically conical, with conic end a homogeneous nearly K\"ahler $6$-manifold. In the case at hand, the end is given by $S^3\times S^3$, but not with the bi-invariant metric.
In this section we use the asymptotically conical model to study the asymptotic behaviour of the connection $A=fA_1$ constructed in the previous section. 

We first note that $S^3$ can be expressed as a symmetric space by the quotient $S^3=(S^3\times S^3)/\Delta$, where the subgroup $\Delta$ is the diagonal subgroup $\Delta\subseteq S^3\times S^3$ acting on the right. This allows us to explicitly identify $S^3\times S^3$ as being the spin bundle $\tilde{\mathcal{F}}$. The fundamental $1$-forms $\omega$ and $\phi$ on $\tilde{\mathcal{F}}$ are then given by $\omega=\frac{1}{2}(\theta_1-\theta_2)$ and $\phi=\frac{1}{2}(\theta_1+\theta_2)$ where the $\theta_i$ are the Maurer-Cartan forms (with values in $\mathfrak{su}(2)$) for the two factors, pulled back to $S^3\times S^3$. Then, the manifold $X$ that supports the $G_2$-metric is given by $X=(S^3\times S^3\times \Ham)/S^3$, with $S^3$ acting again on the right. We will take polar coordinates on $\Ham$, so that the form $\alpha =da-a\phi$ determines the symmetric tensor 
\begin{eqnarray*}
(\alpha^0)^2+(\alpha^1)^2+(\alpha^2)^2+(\alpha^3)^2&=&dr^2+r^2(\theta_3-\phi)^2
\end{eqnarray*}
where $\theta_3$ denotes the Maurer-Cartan form on $S^3\subseteq \Ham$, that being the third factor in the expression $S^3\times S^3 =(S^3\times S^3\times S^3)/S^3$.

The metric $g_\gamma$ of Bryant and Salamon is then given by
\begin{eqnarray*}
g_\gamma&=& 3(1+r^2)^{2/3}\left( (\omega^1)^2+(\omega^2)^2+(\omega^3)^2\right) +4(1+r^2)^{-1/3}\left((\alpha^0)^2+(\alpha^1)^2+(\alpha^2)^2+(\alpha^3)^2\right)\\
&=& 3(1+r^2)^{2/3}\omega^2 +4(1+r^2)^{-1/3}\left(dr^2+r^2(\theta_3-\phi)^2\right).
\end{eqnarray*}
For the squared terms $\omega^2$ and $(\theta_3-\phi)^2$ in this expression, we take sum of the symmetric products of the coefficients of $i$, $j$ and $k$ in each case. Written in this form, it becomes evident that the metric is of cohomegeneity one, with the group $S^3\times S^3\times S^3$ acting on $S^3\times S^3\times \Ham$ on the left and preserving the level sets of $r$. These transformations preserve the forms $\theta_i$, and hence also $\phi$. 
If we use the change of radial variable $\rho=3(1+r^2)^{1/3}$ one sees that 
\begin{eqnarray*}
4(1+r^2)^{-1/3}dr^2&=&\frac{1}{1-(\frac{3}{\rho})^3}d\rho^2,\\
4(1+r^2)^{-1/3}r^2&=& \frac{4}{9}\rho^2\big(1-(\frac{3}{\rho})^3\big),\\
3(1+r^2)^{2/3}&=& \frac{1}{3}\rho^2
\end{eqnarray*}
and so the metric can be expressed as
\begin{eqnarray*}
g_\gamma&=& \frac{1}{1-(\frac{3}{\rho})^3}d\rho^2+\frac{4}{9}\rho^2\left(1-\left(\frac{3}{\rho}\right)^3\right) \left(\theta_3-\phi\right)^2+\frac{1}{3}\rho^2\omega^2.
\end{eqnarray*}
In this form we can recognise the asymptotic behaviour of the metric. As $\rho$ goes to infinity, the metric $g_\gamma$ converges to a conical metric 
\begin{eqnarray*}
g_{con}&=& d\rho^2+\rho^2\Big(\frac{4}{9}(\theta_3-\phi)^2+\frac{1}{3}\omega^2\Big).
\end{eqnarray*}
This is to say that the metric of Bryant and Salamon is asymptotically conical, with the limit being a cone over $S^3\times S^3$ equipped with the metric
\begin{eqnarray}\label{eqn:metriconend}
g=\frac{4}{9}(\theta_3-\phi)^2+\frac{1}{3}\omega^2.
\end{eqnarray}
One can furthermore say that the rate of convergence to the conical model is $O(\rho^{-3})$ \cite{KL}. 

\begin{prop} 
Let $(M,g)$ be an oriented $6$-dimensional Riemannian manifold.  A reduction of the structure group from $SO(6)$ to $SU(3)$ is determined by forms $\varpi\in\Omega^2(M)$, $\Omega_1,\Omega_2\in\Omega^3(M)$ that satisfy the relations 
\begin{eqnarray*}
\frac{\varpi^3}{3!}&=& \frac{1}{4}\Omega_1\wedge\Omega_2= \text{d}vol_g\\
\varpi\wedge\Omega_1&=&\varpi\wedge\Omega_2 =0.
\end{eqnarray*}
\end{prop}
\proof{See for example, \cite{H}}.

This is to say that the forms $\Omega_1$ and $\Omega_2$ determine an almost-complex structure on $M$ compatible with the metric such that $\Omega=\Omega_1+i\Omega_2$ is of type $(3,0)$ and $\varpi$ is the associated Hermitian form.  

\begin{defn} Let $(M,g)$ be a $6$-dimensional Riemannian manifold. The Riemannian cone $C(M)$ over $(M,g)$ is given by the metric $h=dr^2+r^2g$ on $\Real^+\times M$. 
A \emph{conical $G_2$-metric} is a Riemannian cone metric with holonomy contained in $G_2$. 
\end{defn}
\begin{prop}\label{prop:conicalforms}
Let $(C(M),h)$ be a conical $G_2$-metric. The fundamental $3$ and $4$ forms can be expressed as 
\begin{eqnarray*}
\varphi&=&r^2dr\wedge\varpi +r^3\Omega_1\\
*\varphi &=& -r^3dr\wedge \Omega_2+ r^4\frac{\varpi^2}{2}
\end{eqnarray*}
where the differential forms $\varpi\in\Omega^2(M)$, $\Omega_1,\Omega_2\in\Omega^3(M)$ define an $SU(3)$ structure on $M$ and satisfy
\begin{eqnarray}\label{eqn:nearlykahler}
d\Omega_2=-2\varpi^2,\ \ \ d\varpi=3\Omega_1.
\end{eqnarray}
\end{prop}
\proof{ The forms are defined by contracting by the radial vector field, whose flow is by homotheties. For $\varphi$ to be closed and coclosed, the forms must be related as in Equation (\ref{eqn:nearlykahler}). }
\begin{defn}
A ($6$-dimensional) \emph{Nearly-K\"ahler manifold} is a Riemannian manifold $(M^6,g)$ with structure group reduced to $SU(3)$ by forms $\varpi,\Omega_1,\Omega_2$ that satisfy Equation (\ref{eqn:nearlykahler}).
\end{defn}
\begin{prop} The three and four-forms $\gamma$ and $*\gamma$ on $S^3\times \Ham$ given by the Bryant-Salamon construction are expressed in the new coordinate system as 
\begin{eqnarray}
\gamma &=& \frac{1}{3\sqrt{3}}\rho^3\omega^{123} +2\sqrt{3}\Big(\left(\frac{\rho}{3}\right)^3-1\Big) \text{Re}\Big(\omega\wedge(\theta_3-\phi)^2\Big) \\
&&\ \ \ \ \ \ \ \ \ \ \ \ \ \ \ \ \ \ \ +\frac{2}{3\sqrt{3}}\rho^2d\rho\wedge \text{Re}\Big(\omega\wedge (\theta_3-\phi)\Big)\label{eqn:gammaasymp}\nonumber\\
*\gamma &=& \frac{-4}{3}\Big(\left(\frac{\rho}{3}\right)^3-1\Big)d\rho\wedge (\theta_3- \phi)^3 +\frac{1}{9}\rho^3d\rho\wedge \text{Re}\Big(\omega^2\wedge (\theta_3-\phi)\Big)\\
&& \ \ \ \ \ \ \ \ \ \ \ \ \ \ \ \ \ \ -\rho\Big(\left(\frac{\rho}{3}\right)^3-1\Big) \text{Re}\Big(\omega^2\wedge(\theta_3-\phi)^2\Big). \nonumber
\end{eqnarray}
As $\rho$ goes to infinity, these forms respectively converge to 
\begin{eqnarray*}
\gamma_{con} &=& \rho^2 d\rho\wedge \left[ \frac{2}{3\sqrt{3}} \text{Re} \Big(\omega\wedge(\theta_3-\phi)\Big)\right] + \rho^3 \left[ \frac{1}{3\sqrt{3}}\omega^{123} +\frac{2}{9\sqrt{3}} \text{Re} \Big(\omega\wedge (\theta_3-\phi)^2\Big)\right]\\
*\gamma_{con} &=& \rho^3 d\rho\wedge \left[\frac{-4}{81}(\theta_3-\phi)^3+\frac{1}{9}\text{Re}\Big(\omega^2\wedge(\theta_3-\phi)\Big)\right] -\rho^4\left[\frac{1}{27}\text{Re}\Big(\omega^2\wedge(\theta_3-\phi)^2                             \Big)\right].
\end{eqnarray*}
\end{prop}
We note here that in the squared terms here, for example $\omega^2$, we are taking the exterior product of $\text{Im}(\Ham)$-valued $1$-forms, rather than the symmetric product as above.
\proof{In the following calculations we repeatedly use the identity $\bar a\alpha=\bar ada-r^2\phi= d(\frac{r^2}{2})+r^2(\theta_3-\phi)$, which can be readily verified. We recall that $\gamma=f^3\gamma_1+fg^2\gamma_2$ where $\gamma_1=\omega^{123}$ and $\gamma_2=\frac{-1}{2}\text{Re}(\omega\wedge\bar\alpha\wedge\alpha)$. We also observe that $fg^2=4\sqrt{3}$ and $f^3=3\sqrt{3}(1+r^2)=\frac{1}{3\sqrt{3}}\rho^3$. Then, with $r^2=(\rho/3)^3-1$ and $d(r^2)=\frac{\rho^2}{9}d\rho$,
\begin{eqnarray*}
r^2\gamma_2&=& \frac{-1}{2}\text{Re}\left(\omega\wedge\bar\alpha a\wedge\bar a\alpha\right)\\
&=& \frac{-1}{2}\text{Re}\left(\omega\wedge\Big(d(\frac{r^2}{2})-r^2(\theta_3-\phi)\Big)\wedge \Big(d(\frac{r^2}{2})+r^2(\theta_3-\phi)\Big)\right)\\
&=& \frac{r^2}{2}d(r^2)\wedge\text{Re}\left(\omega\wedge(\theta_3-\phi)\right)+\frac{r^4}{2}\text{Re}\left(\omega\wedge(\theta_3-\phi)^2\right)\\
4\sqrt{3}\gamma_2&=& \frac{4\sqrt{3}}{2\cdot 9}\rho^2 d\rho\wedge\text{Re}\left(\omega\wedge(\theta_3-\phi)\right) +2\sqrt{3}\left(\Big(\frac{\rho}{3}\Big)^3-1\right)\text{Re}\left(\omega\wedge(\theta_3-\phi)^2\right)
\end{eqnarray*}
which is sufficient to prove Equation (\ref{eqn:gammaasymp}). The asymptotic limit is obtained by taking the dominant powers of $\rho$, and observing that the remainder tends to zero in the metric $g_{con}$ for large $\rho$.

For the other term $\psi=*\gamma$, we first recall that $g^4=16(1+r^2)^{-2/3}=16\cdot 9/\rho^2$ and that $f^2g^2=3.4(1+r^2)^{2/3}(1+r^2)^{-1/3}=4\rho$. Then,
\begin{eqnarray*}
\psi_1&=& \alpha^{0123} =\frac{1}{24}\bar\alpha\wedge\alpha\wedge\bar\alpha\wedge\alpha\\
24r^4\alpha^{0123}&=& \left[\Big(d(\frac{r^2}{2})-r^2(\theta_3-\phi)\Big)\wedge\Big(d(\frac{r^2}{2})+r^2(\theta_3-\phi)\Big)\right]^2\\
&=& -2r^6d(r^2)\wedge(\theta_3-\phi)^3\\
\alpha^{0123} &=& \frac{-r^2}{12}d(r^2)\wedge(\theta_3-\phi)^3\\
&=& \frac{-1}{12}\frac{\rho^2}{9}\Big((\frac{\rho}{3})^3-1\Big)d\rho\wedge(\theta_3-\phi)^3\\
g^4\psi_1&=& \frac{-4}{3}\Big((\frac{\rho}{3})^3-1\Big)d\rho\wedge(\theta_3-\phi)^3.
\end{eqnarray*}
\begin{eqnarray*}
\psi_2&=& \frac{-1}{4}\text{Re}\left(\omega^2\wedge\bar \alpha\wedge\alpha\right)\\
r^2\psi_2&=& \frac{-1}{4}\text{Re}\left(\omega^2\wedge\bar\alpha a\wedge\bar a\alpha\right)\\
&=& \frac{-1}{4}\text{Re}\left(\omega^2\wedge 2r^2d(\frac{r^2}{2})\wedge(\theta_3-\phi)-\omega^2\wedge r^4(\theta_3-\phi)^2\right)\\
\psi_2&=& \frac{-1}{4}(\frac{\rho}{3})^2d\rho\wedge\text{Re}\left(\omega^2\wedge(\theta_3-\phi)\right) +\frac{1}{4}\left((\frac{\rho}{3})^3-1\right)\text{Re}\left(\omega^2\wedge(\theta_3-\phi)^2\right)\\
f^2g^2\psi_2 &=&\frac{-1}{9}\rho^3d\rho\wedge\text{Re}\Big(\omega^2\wedge(\theta_3-\phi)\Big)+\rho\Big((\frac{\rho}{3})^3-1\Big)\text{Re}\left(\omega^2\wedge(\theta_3-\phi)^2\right).
\end{eqnarray*}
Putting these terms together as $\psi=g^4\psi_1-f^2g^2\psi_2$ gives the desired formula for $\psi=*\gamma$. 
}

As such, we can define the forms $\varpi$, $\Omega_1$ and $\Omega_2$ that define an $SU(3)$ structure on $S^3\times S^3$ :
\begin{eqnarray}
\varpi &=&  \frac{2}{3\sqrt{3}}\text{Re}\Big(\omega\wedge(\theta_3-\phi)\Big)\nonumber\\
\Omega_1 &=&\frac{1}{3\sqrt{3}}\omega^{123} +\frac{2}{9\sqrt{3}}\text{Re}\Big(\omega\wedge (\theta_3-\phi)^2\Big)\nonumber\\
\Omega_2 &=& \frac{-4}{81}(\theta_3-\phi)^3+\frac{1}{9}\text{Re}\Big(\omega^2\wedge(\theta_3-\phi)\Big)\label{eqn:holomvolform2}
\end{eqnarray}
Furthermore one can show without trouble that 
\begin{eqnarray}\label{eqn:omegasquared}
\frac{1}{2}\varpi^2 &=& \frac{-1}{27}\text{Re}\Big(\omega^2\wedge(\theta_3-\phi)^2  \Big)
\end{eqnarray}
as dictated by the decomposition of the conical $4$-form in Proposition \ref{prop:conicalforms}. If we return to the definition of the metric on $S^3\times S^3$ in Equation (\ref{eqn:metriconend}) we can give an orthonormal coframe. If we write $\omega=\omega^1i+\omega^2j+\omega^3k$ and $\theta_3-\phi=\tau^1i+\tau^2j+\tau^3k$, then setting $\sigma^\beta=1/\sqrt{3}\omega^\beta$ (for $\beta=1,2,3$) and $\sigma^\beta=2/3\tau^{\beta-3}$ (for $\beta=4,5,6$) gives an orthonormal coframe on $S^3\times S^3$. The volume form is then given by $dvol=\frac{-1}{3\sqrt{3}}(\frac{2}{3})^3\omega^{123}\wedge\tau^{123}$. The differential forms $\Omega_1$ and $\Omega_2$ are then given by $\Omega_1=\text{Re}(\Omega)$ and $\Omega_2=\text{Im}(\Omega)$ where 
\begin{eqnarray*}
\Omega=\eta^1\wedge\eta^2\wedge\eta^3
\end{eqnarray*}
for $\eta^1=\sigma^1+i\sigma^4$, $\eta^2=\sigma^2+i\sigma^5$ and $\eta^3=\sigma^3+i\sigma^6$.

We now return to the connection $A=f(r)A_1$ constructed on a trivial $SU(2)$-bundle over $X$. 
\begin{defn} 
Let $(M^6,g)$ be a Nearly K\"ahler manifold with structure determined by the forms $\varpi$, $\Omega=\Omega_1+i\Omega_2$. Let $\mathcal{P}\to M$ be a principal $G$-bundle on $M$. A connection $A$ on $\mathcal{P}$ is \emph{Hermitian Yang-Mills} if its curvature satisfies the conditions
\begin{eqnarray*}
F_A\wedge \varpi^2 &=& 0\\
F_A\wedge \Omega&=& 0.
\end{eqnarray*}
\end{defn}
This is to say that the curvature is a primitive $(1,1)$-form on $M$, with respect to the almost-complex structure and hermitian form given on $M$. Hermitian Yang Mills connections were first defined and studied on non-integrable almost-complex manifolds by Bryant \cite{B}, and subsequently studied on Nearly K\"ahler manifolds by Xu \cite{X} and Oliveira \cite{O}.

$A_1$ is given by $A_1=\text{Im}(a\bar \alpha)$ so
\begin{eqnarray*}
r^2A_1=\text{Im}(a\bar\alpha a\bar a)&=&\text{Im}(a(d(\frac{r^2}{2})-r^2(\theta_3-\phi))\bar a)\\
&=&-r^2a(\theta_3-\phi)\bar a\\
&=& -r^4 g_3(\theta_3-\phi)g_3^{-1}
\end{eqnarray*}
for $(g_1,g_2,g_3)\in S^3\times S^3\times S^3$ so $A_1=-r^2 g_3(\theta_3-\phi)g_3^{-1}=-((\frac{\rho}{3})^3-1)g_3(\theta_3-\phi)g_3^{-1}$. This form is defined on $S^3\times S^3\times S^3$ but is invariant by the right action of the diagonal subgroup $\Delta$ and vanishes along the orbits of that group. We will in the following write $g$ for $g_3$. The connection $A$ constructed in the previous section is then given by
\begin{eqnarray*}
A=f(r)A_1 =\frac{-2r^2}{3(1+r^2)+C(1+r^2)^{1/3}}g(\theta_3-\phi)g^{-1}=\frac{-2((\frac{\rho}{3})^3-1)}{\frac{1}{9}\rho^3+\frac{C}{3}\rho}g(\theta_3-\phi)g^{-1}.
\end{eqnarray*}
This connection, defined as a differential form on $S^3\times S^3\times \Ham$, is invariant under the right action of the diagonal group $\Delta\subseteq S^3\times S^3\times S^3$, but it is also invariant under the left action of the full group $(S^3)^3$, and so has the same level of symmetry (or at least the same identity component) as the metric $g_\gamma$.
As $\rho$ goes to infinity, this connection converges on the conical end to the pull-back from the slice $S^3\times S^3$ of the connection
\begin{eqnarray*}
\widetilde{A}=\frac{-2}{3}g(\theta_3-\phi)g^{-1}.
\end{eqnarray*}
We can now restate and prove Theorem \ref{thm:asymptotics}.
\begin{thm}
The asymptotic limit connection $\widetilde{A}$ is a Hermitian Yang-Mills connection on the trivial $SU(2)$-bundle over $S^3\times S^3$.
\end{thm}
\proof{ We can use the identities $dg=g\theta_3$, $dg^{-1}=-\theta_3 g^{-1}$, $d\theta_3=-\theta_3\wedge\theta_3$ and $d\phi=-\phi\wedge\phi-\frac{1}{4}\omega\wedge\omega$ to show that the curvature form $\widetilde{A}$ satisfies 
\begin{eqnarray*}
d\widetilde{A} &=& \frac{-6}{9}g(\theta_3-\phi)^2g^{-1}-\frac{1}{6}g\omega^2 g^{-1},\\
\widetilde{A}\wedge\widetilde{A} &=& \frac{4}{9}g(\theta_3-\phi)^2g^{-1},\\
F_{\widetilde{A}}&=& \frac{-2}{9}g(\theta_3-\phi)^2g^{-1}-\frac{1}{6}g\omega^2 g^{-1}.
\end{eqnarray*}
Firstly, from Equation (\ref{eqn:omegasquared}) it is clear that $F_{\widetilde{A}}\wedge\omega^2=0$. For $\Omega_2$ from Equation (\ref{eqn:holomvolform2}), we  then see that
\begin{eqnarray*}
9\Omega_2 &=& \frac{8}{3}\tau^{123}-2\left(\omega^{23}\wedge\tau^1+\omega^{31}\wedge \tau^2+\omega^{12}\wedge\tau^3\right)\\
9F_{\widetilde{A}}\wedge\Omega_2 &=& \frac{-8}{18}\tau^{123}\wedge\omega^2 +\frac{8}{9}(\omega^{23}i+\omega^{31}j+\omega^{12}k)\wedge\tau^{123}\\
&=& \frac{-8}{18}\tau^{123}\wedge\omega^2+\frac{8}{9}\frac{1}{2}\omega^2\wedge\tau^{123}\\
&=&0
\end{eqnarray*}
with a similar calculation for $\Omega_1$. This completes the proof of Theorem \ref{thm:asymptotics}.
}

\section{$Spin(7)$-Instantons on the negative spinor bundle of $S^4$} \label{sec:Spin7inst}

In this section we prove Theorem \ref{thm:Spin7} and construct explicit $Spin(7)$-instantons with structure group $SU(2)$ on the $8$-dimensional manifold of Bryant and Salamon
 equipped with the $Spin(7)$-holonomy metric. That is, let $\mathcal{F}$ be the spin structure for the round metric on $S^4$, and $Y=(\mathcal{F}\times \Ham)/Spin(4)$ the negative spinor bundle, considered as an $8$-dimensional manifold. $\mathcal{P}=\mathcal{F}\times \Ham$ is a principal $Spin(4)=SU(2)\times SU(2)$ bundle and $\mathcal{E}_Y=(\mathcal{F}\times \Ham)/(SU(2)\times \{1\})$ is a principal $SU(2)$-bundle over $Y$. 

The Levi-Civita connection $\theta=\psi+\phi$ pulls back to $\mathcal{F}\times \Ham$ and the form $\phi$ descends to $\mathcal{E}_Y$ to define a connection with structure group $SU(2)$.
 Let $A_1$ be the $\mathfrak{su}(2)=\text{Im}(\Ham)$-valued $1$-form on $\mathcal{P}$ defined by 
\begin{eqnarray*}
A_1=\text{Im}(\bar a\alpha)=\frac{1}{2}(\bar a\alpha-\bar \alpha a)
\end{eqnarray*}
where $a$ and $\alpha$ are defined in Section \ref{Sect:BSmetric}. $A_1$ descends to an equivariant form on $\mathcal{E}_Y$ but all calculations will be performed on $\mathcal{P}$. Then, $A_1$ is equal to zero on vectors tangent to the fibres of $\mathcal{P}\to X$ and is equivariant in the sense that 
\begin{eqnarray*}
R_g^*A_1 =\bar qA_1q
\end{eqnarray*}
for $g=(p,q)\in SU(2)\times SU(2)$. 

For any $SU(2)$-invariant function on $\mathcal{E}_Y$, the form $A=\phi+fA_1$ defines a connection on $\mathcal{E}_Y$. We will suppose that $f$ depends only on $r=\bar a a$. 
The curvature of $A$ is equal to 
\begin{eqnarray*}
F_A&=& \left[d\phi +\phi\wedge\phi\right] +f^\prime dr\wedge A_1+f\left[dA_1+A_1\phi+\phi A_1\right] +f^2A_1\wedge A_1.
\end{eqnarray*}

\begin{lemma} We can simplify the forms as follows.
\begin{eqnarray}
dA_1+A_1\phi+\phi A_1&=& \bar \alpha\wedge\alpha -\frac{r\kappa}{2}\Omega,\\
A_1\wedge A_1&=& \frac{-r}{2}\bar\alpha\wedge\alpha -\bar a(\frac{1}{2}\alpha\wedge\bar\alpha)a,\\
dr\wedge A_1 &=& \frac{r}{2}\bar\alpha\wedge\alpha - \bar a(\frac{1}{2}\alpha\wedge\bar\alpha)a.
\end{eqnarray}
\end{lemma}

\proof{
We first recall the formulae from Bryant and Salamon, 
\begin{eqnarray*}
d\alpha+\alpha\phi&=& -a\frac{\kappa}{2}\Omega\\
d\phi+\phi\wedge\phi &=& \frac{\kappa}{2}\Omega
\end{eqnarray*}
where $\Omega=1/2\bar\omega\wedge\omega$. Then,
\begin{eqnarray*}
2dA_1 &=& (\bar\alpha-\phi\bar a)\alpha +\bar a(-\alpha\phi-\frac{\kappa}{2}a\Omega)-(\frac{\kappa}{2}\Omega\bar a-\phi\bar\alpha)a +\bar \alpha(\alpha+ a\phi)\\
&=& 2\bar\alpha\wedge\alpha -\phi(\bar a\alpha-\bar \alpha a)-(\bar a\alpha-\bar\alpha a)\phi-r\kappa\Omega.
\end{eqnarray*}
For the second identity,
\begin{eqnarray*}
4A_1\wedge A_1 &=& (\bar a\alpha-\bar\alpha a)\wedge(\bar a\alpha-\bar\alpha a)\\
&=& \bar a\alpha \bar a\alpha-\bar \alpha \bar aa\alpha-\bar a\alpha\bar\alpha a+\bar \alpha a\bar\alpha a\\
&=& 2A_1\wedge A_1 -r\bar \alpha\wedge\alpha-\bar a\alpha\wedge\bar\alpha a.
\end{eqnarray*}
This calculation follows from the observation that if $B=\text{Im}(\varphi)$ then $B\wedge B=\varphi\wedge\varphi=\bar\varphi\wedge\bar\varphi$ if $\varphi$ is a quaternion-valued $1$-form. Finally,
\begin{eqnarray*}
2dr\wedge A_1&=&(\bar a\alpha +\alpha\bar a)\wedge(\bar a\alpha -\bar\alpha a)\\
&=& \bar a\alpha\bar a\alpha+r\bar\alpha\alpha -\bar a\alpha\bar\alpha a-\bar\alpha a\bar\alpha a\\
dr\wedge A_1&=& \frac{r}{2}\bar\alpha\wedge\alpha -\bar a(\frac{1}{2}\alpha\wedge\bar\alpha)a.
\end{eqnarray*}
}

With respect to this decomposition, the curvature term $F_A$ can be expressed as 
\begin{eqnarray*}
F_A &=& \frac{\kappa}{2}\left(1-rf\right)\Omega +\left(rf^\prime +2f-rf^2\right)\frac{1}{2}\bar\alpha\wedge\alpha \\ 
&&\ \ \ \  \ -\left(f^\prime +f^2\right)\bar a(\frac{1}{2}\alpha\wedge\bar\alpha)a.
\end{eqnarray*}
We denote the three terms by 
\begin{eqnarray*}
F_1 &=& \Omega=\Omega^1i+\Omega^2j+k\Omega^3\\
F_2 &=& \frac{1}{2}\bar\alpha\wedge\alpha =B= B^1i+B^2j+B^3k\\
F_3 &=& \bar a(\frac{1}{2}\alpha\wedge\bar\alpha)a
\end{eqnarray*}
A connection $A=\phi+fA$ is a $Spin(7)$-instanton if $\Psi\wedge F_A=*F_A$ where $*$ is the Hodge dual operator. The fundamental $4$-form on $X$ is given by $\Psi=\sigma^2\Psi_1+\sigma\tau\Psi_2+\tau^2\Psi_3$ where
\begin{eqnarray*}
\Psi_1&=&\frac{-1}{6}B\wedge B=\alpha^{0123},\\
\Psi_2 &=&-\text{Re}(B\wedge\Omega)= B^1\wedge\Omega^1+B^2\wedge\Omega^2+B^3\wedge\Omega^3,\\
\Psi_3 &=& \frac{-1}{6}\Omega\wedge\Omega =\omega^{0123}.
\end{eqnarray*}
and where $\sigma(r)=4(1+r)^{-2/5}$ and $\tau(r)=5\kappa(1+r)^{3/5}$. 

\begin{lemma}
The products satisfy $\Psi_1\wedge F_2=\Psi_1\wedge F_3=\Psi_2\wedge F_3=\Psi_3\wedge F_1=0$. The other terms satisfy
\begin{eqnarray*}
\Psi_1\wedge F_1&=&\Omega\wedge\alpha^{0123},\\
\Psi_2\wedge F_2&=& -2\Omega\wedge \alpha^{0123}\\
\Psi_2\wedge F_1 &=& -2B\wedge \omega^{0123}\\
\Psi_3\wedge F_2 &=& B\wedge \omega^{0123}\\
\Psi_3\wedge F_3 &=& \bar a (\frac{1}{2}\alpha\wedge\bar\alpha) a\wedge\omega^{0123}.
\end{eqnarray*}
\end{lemma}
Apart from the cases that are immediate, this follows from the observation that $B=B^1i+B^2j+B^3k$ where $(B^i)^2=-2\alpha^{0123}$, and similarly for the $\Omega$ terms. It is evident that the three $6$-forms $B\wedge \omega^{0123}$, $\Omega\wedge\alpha^{0123}$ and $\bar a\alpha\wedge\bar\alpha a\wedge\omega^{0123}$ are linearly independent at each point because, for example, $B$ has ``anti-self-dual" components, while $\bar a\alpha\wedge\bar\alpha a$ has ``self-dual" components.

The above calculations demonstrate the following proposition.
\begin{prop}
\begin{eqnarray}
\Psi\wedge F_A&=& 
\Big( \sigma^2\frac{\kappa}{2}(1-rf) -2\sigma\tau (rf^\prime +2f-rf^2)\Big)\Omega\wedge\alpha^{0123}\label{eqn:formwedgecurve}\\ 
&& \ \ \ +\ \Big(-2\sigma\tau\frac{\kappa}{2}(1-rf) +\tau^2(rf^\prime +2f-rf^2)\Big)B\wedge\omega^{0123}\nonumber\\
&& \ \ \ \ \ \ \ \ \ -\ \tau^2\Big(f^\prime +f^2\Big)\bar a(\frac{1}{2}\alpha\wedge\bar\alpha)a\wedge\omega^{0123}.\nonumber
\end{eqnarray}
\end{prop}

The metric that is determined by $\Psi$ is given by
\begin{eqnarray*}
g_\Psi = \sigma\left((\alpha^0)^2+(\alpha^1)^2+(\alpha^2)^2+(\alpha^3)^2\right)+\tau\left((\omega^0)^2+ (\omega^1)^2+(\omega^2)^2+(\omega^3)^2\right)
\end{eqnarray*}
and so an orthonormal coframe is given by $\{\sqrt{\sigma}\alpha^i,\sqrt{\tau}\omega^j\}$. The Hodge star operator therefore satisfies 
\begin{eqnarray*}
*(\bar\omega\wedge\omega)&=& -\sigma^2\bar\omega\wedge\omega\wedge\alpha^{0123},\\
*(\bar\alpha\wedge\alpha    )&=& - \tau^2\bar\alpha\wedge\alpha\wedge\omega^{0123},\\
*(\bar a\alpha\wedge\bar\alpha a) &=& \tau^2\bar a\alpha\wedge\bar\alpha a\wedge\omega^{0123}.
\end{eqnarray*}

\begin{prop}
\begin{eqnarray}
*F_A &=& -\frac{\sigma^2\kappa}{2}(1-fr)\Omega\wedge\alpha^{0123} \label{eqn:starcurve}\\
&& \ \ \ \ \ \ \ \ \ -\ \tau^2\Big(rf^\prime +2f-rf^2\Big)B\wedge\omega^{0123}\nonumber\\
&&\ \ \ \ \ \ \ \ \ \ \ \ \ \ -\ \tau^2\Big(f^\prime+f^2\Big)\bar a(\frac{1}{2}\alpha\wedge\bar\alpha )a\wedge \omega^{0123}.\nonumber
\end{eqnarray}
\end{prop}

We can note that for any function $f$, the terms involving $\bar a(\frac{1}{2}\alpha\wedge\bar\alpha)a$ in Equations \ref{eqn:formwedgecurve} and \ref{eqn:starcurve} coincide. We can therefore see that $A$ defines a $Spin(7)$-instanton if the other terms coincide as well. By pure luck, perhaps, these two conditions give the same equation. That is, 
\begin{eqnarray*}
\frac{\sigma^2\kappa}{2}(1-rf)-2\sigma\tau\Big(rf^\prime+2f-rf^2\Big)=\frac{-\sigma^2\kappa}{2}(1-rf).
\end{eqnarray*}

These two facts; the coincidence of the coefficients of $\bar a\frac{\alpha\wedge\bar\alpha}{2} a$ and the two conditions reducing to one equation; suggest a great deal of symmetry in this system of equations. For any functions $\sigma$ and $\tau$ the metric constructed on $Y$ is of cohomogeneity-one, and the isometries lift to automorphisms of the bundle $\mathcal{E}_Y$ that preserve $A$. The author plans to make more precise this statement, and to study these connections on the asymptotically conical model for $Y$, as in the $G_2$ case earlier.

Then, using the above choices for $\sigma$ and $\tau$, we can conclude the following proposition.

\begin{prop}
The connection $A$ satisfies $\Psi\wedge F_A=*F_A$ if and only if $f$ satisfies
\begin{eqnarray}\label{eqn:fdiffeqn}
rf^\prime +\frac{12r+10}{5(1+r)} f-rf^2= \frac{2}{5}\frac{1}{1+r}.
\end{eqnarray}
\end{prop}

As in the previous example, for an instanton on the $G_2$-manifold constructed by Bryant and Salamon, we can see that this ordinary differential equation is of Riccati type. The principal difference in this case is that there exists an inhomogeneous term, that comes in this case from the fact that the base point connection $\phi$ is not flat, as it was in the previous case. This can be overcome as follows. 

\begin{prop}
Let $f$ satisfy Equation \ref{eqn:fdiffeqn}. Then the function $g=f-\frac{1}{r}$ satisfies 
\begin{eqnarray}\label{eqn:gdiffeqn}
g^\prime +\frac{2}{5(1+r)}g-g^2=0.
\end{eqnarray}
\end{prop}

This substitution can be motivated geometrically. In establishing the ansatz above, we took the base connection on $\mathcal{E}_Y$ to be the form $\phi$ pulled back from $\mathcal{F}$. We can instead take the form $\varphi=\text{Im}(a^{-1}da)$, which is defined away from the zero section $\{a=0\}$. That is, $\varphi=\phi+1/r\text{Im}(\bar a\alpha)$, and so $\varphi +gA_1$ is a $Spin(7)$-instanton if $g$ satisfies Equation \ref{eqn:gdiffeqn}. We can first note that $g=0$ gives a solution, which shows that $\varphi=\text{Im}(a^{-1}da)$ defines an instanton away from $\{r=0\}$. 

As noted above, Equation \ref{eqn:gdiffeqn} is of Riccati type, so we can make a similar substitution. Suppose that $g(r)=-y^\prime(r)/y(r)$. Then $y$ satisfies the equation
\begin{eqnarray*}
y^{\prime\prime}+\frac{2}{5}\frac{1}{1+r}y^\prime&=&0.
\end{eqnarray*}
This yields solutions
\begin{eqnarray*}
y(r)&=&A(1+r)^{3/5} +B,\\
g(r)&=& \frac{-3C}{(1+r)^{2/5}(1+5C(1+r)^{3/5})},\\
f(r)&=& \frac{1}{r(1+D(1+r)^{3/5})}+\frac{D(2r+5)}{5r(1+r)^{2/5}(1+D(1+r)^{3/5})}.
\end{eqnarray*}
If $D$ is taken to be greater than $-1$, this function is singular only for $r=0$, and the connection $A$ is defined everywhere away from the zero section locus of $Y$. We note that $D=0$ corresponds to the function $f=1/r$ or $g=0$. This completes the proof of Theorem \ref{thm:Spin7}.

\subsection*{Acknowledgements} The author would like to thank the reviewer for many helpful suggestions that greatly improved the paper. He would also like to acknowledge the financial support of FAPESP processo 2011/07363-6.

\bibliographystyle{ams}

\begin{thebibliography}{EFK}

\bibitem[AOS]{AOS} B.S. Acharya, M. O'Loughlin and B. Spence, {\em Higher dimensional analogues of Donaldson-Witten theory}, Nuclear Phys. B {\bf 503} (1997), no. 3, 657-674.

\bibitem[AW]{AW} M. Atiyah and E. Witten, {\em $M$-theory dynamics on a manifold of $G_2$ holonomy}, Adv. Theor. Math. Phys., {\bf 6} (2002), 1-106.

\bibitem[BKS]{BKS} L. Baulieu, H. Kanno and I.M. Singer, {\em Special quantum field theories in eight and other dimensions}, Comm. Math. Phys. {\bf 194} (1998) no. 1, 149-175.

\bibitem[BPST]{BPST} A.A. Belavin, A.M. Polyakov, A.S. Schwartz and Yu.S. Tyurin, {\em Pseudo-particles solutions of the Yang-Mills equations}, Phys. Lett. B., {\bf 59}, 87--89.


\bibitem[B]{B} R.L. Bryant, {\em On the geometry of almost-complex $6$-manifolds}, Asian J. Math., Vol. 10, No. 3 (2006), 561-606.

\bibitem[BS]{BS} R.L. Bryant and S. M. Salamon, {\em On the construction of some complete metrics with exceptional holonomy}, Duke Math. J.,  {\bf 58} (1989), 829-850.


\bibitem[CGLP]{CGLP} M. Cvetic, G.W. Gibbons, H. Lu and C.N. Pope, {\em Cohomogeneity one manifolds of $Spin(7)$ and $G_2$ holonomy}, Phys. Rev. D, Vol. 65, No. 10, (2002).

\bibitem[DT]{DT} S.K. Donaldson and R.P. Thomas, {\em Gauge theory in higher dimensions}, in : The Geometric Universe, Oxford Univ. Press, 1998, pp 31-47. 

\bibitem[FN]{FN} D.B. Fairlie and J. Nuyts, {\em Spherically symmetric solutions of gauge theories in eight dimensions}, J. Phys. A {\bf 17} (1984), no. 14, 2867-2872.

\bibitem[GPP]{GPP} G.W. Gibbons, D.N. Page and C.N. Pope, {\em Einstein metrics on $S^3$, $R^3$ and $R^4$ bundles}, Comm. Math. Phys. {\bf 127} (1990), no. 3, 529-553.

\bibitem[HN]{HN} D. Harland and C. N\"olle, {Instantons and Killing spinors}, J. High Energy Phys., (2012), No. 3, 082.

\bibitem[HILP]{HILP} D. Harland. T.A. Ivanova, O. Lechtenfeld and A.D. Popov, {\em Yang-Mills flows on nearly K\"aher manifolds and $G_2$-instantons}, Comm. Math. Phys. 300:185-204 (2010).

\bibitem[H]{H} N. Hitchin, {\em The geometry of three-forms in six and seven dimensions}, J. Diff. Geom. {\bf 55} (2000), no. 3, 547--576.

\bibitem[IP]{IP} T.A. Ivanova and A.D. Popov, {\em Instantons on special holonomy manifolds}, Phys.Rev D85, 105012 (2012) 

\bibitem[KL]{KL} S. Karigiannis and J. Lotay, {\em Deformation theory of $G_2$-conifolds},  arXiv:1212.6457.



\bibitem[O]{O} G. Oliveira, {\em Monopoles on the Bryant-Salamon $G_2$ manifolds}, preprint, available at arxiv:13107392v1.

\bibitem[S]{S} H.N. S\'a Earp, { Instantons on $G_2$-manifolds},  PhD thesis, Imperial College London, (2009). 


\bibitem[Ti]{Ti} G. Tian, {\em Gauge theory and calibrated geometry, I}, Ann. of Math. (2) {\bf 151} (2000) no. 1, 193-268.

\bibitem[W]{W} T. Walpuski, {\em $G_2$-instantons on generalised Kummer constructions}, Geometry and Topology, {\bf 17} (2013), 2345--2388.

\bibitem[X]{X} F. Xu, {\em On instantons on nearly K\"ahler $6$-manifolds}, Asian J. Math., Vol. 13, No. 4, (2009), 535-568.









\end{thebibliography}

\end{document}